%% file: geometric_inference_arxiv.tex
\documentclass[pre,nofootinbib]{revtex4}

\usepackage{natbib}
\usepackage{amsmath}
\usepackage{graphicx}
\usepackage{mathrsfs}
\usepackage{amsthm}
\usepackage{amscd}
\usepackage{amssymb}
\usepackage{latexsym}
\usepackage{stmaryrd}
\usepackage{wasysym}
\usepackage{enumitem}
\usepackage{verbatim}
\usepackage{nicefrac}

\newcommand{\hyperboloid}[1]{\mathbb{H}_{#1}}

\newcommand{\ContinuousSpace}{M}
\newcommand{\ContinuousMetric}{\rho}
\newcommand{\NodeDensityLpMetric}[1]{d_{#1}}
\newcommand{\QuasiUniformDensity}{q}

\newcommand{\KDE}{\widehat{f}^{n,h}}
\newcommand{\gKDE}{\widehat{f}^{n,h}}
\newcommand{\gKDEA}{\widehat{f}^{n_1,h_1}}
\newcommand{\gKDEB}{\widehat{f}^{n_2,h_2}}
\newcommand{\gKDEAB}{\widehat{f}^{n_1+n_2,h_{12}}}
\newcommand{\gKDEBootA}{\widehat{f}^{n_1,h_1 *}}
\newcommand{\gKDEBootB}{\widehat{f}^{n_2,h_2 *}}
\newcommand{\EmpiricalMeasure}{\widehat{P}_n}
\newcommand{\HyperbolicGaussianParam}{\sigma}

\newcommand{\Fourier}[1]{\mathcal{F}\left[ #1 \right]}
\newcommand{\Helgason}[1]{\mathcal{H}\left[ #1 \right]}

\newcommand{\ProbModel}{\mathcal{P}}
\newcommand{\EstModel}{\widehat{\mathcal{P}}}

\newcommand{\SO}{\mathbb{SO}}
\newcommand{\SL}{\mathbb{SL}}
\newcommand{\R}{\mathbb{R}} 
\newcommand{\C}{\mathbb{C}} 

\begin{document}

%

%

\title{Geometric Network Comparisons}

\author{Dena Asta}
\email{dasta@andrew.cmu.edu}
\affiliation{Department of Engineering and Public Policy, Carnegie Mellon University, Pittsburgh, PA 15213 USA}
\affiliation{Department of Statistics, Carnegie Mellon University, Pittsburgh, PA 15213 USA}
\author{Cosma Rohilla Shalizi}
\email{cshalizi@cmu.edu}
\affiliation{Department of Statistics, Carnegie Mellon University, Pittsburgh, PA 15213 USA}
\affiliation{Santa Fe Institute, 1399 Hyde Park Road, Santa Fe, NM 87501, USA}
\date{26 October 2014; last \LaTeX 'd \today}

\begin{abstract}
  Network analysis has a crucial need for tools to compare networks and assess
  the significance of differences between networks.  We propose a principled
  statistical approach to network comparison that approximates networks as
  probability distributions on negatively curved manifolds.  We outline the
  theory, as well as implement the approach on simulated networks.
\end{abstract}

\maketitle

\section{Introduction}
\label{sec:intro}
\input{introduction}

\section{Motivation and Background}
\label{sec:motivation}
\input{motivation}

\section{Methodology}
\label{sec:methodology}
\input{methodology}

\section{Simulations}
\label{sec:experiments}
\input{experiments}

\section{Conclusion}
\label{sec:conclusion}

\input{conclusion}

\section*{Acknowledgements}

Our work was supported by NSF grant DMS-1418124 and NSF Graduate Research Fellowship under grant DGE-1252522.  We are grateful for valuable
discussions with Carl Bergstrom, Elizabeth Casman, David Choi, Aaron Clauset, Steve Fienberg, Christopher
Genovese, Dmitri Krioukov, Alessandro Rinaldo, Mitch Small, Andrew Thomas, Larry Wasserman, Christopher
Wiggins, and for feedback from seminar audiences at UCLA's Institute for Pure
and Applied Mathematics and CMU's machine learning and social science seminar.

\bibliography{extra,locusts}

\clearpage

\appendix

\section{Helgason Transforms}
\label{app:helgason}

\input{helgason}

\section{Efficient Computation of the Test Statistic}
\label{app:test-static-computing}
\input{computing}

\end{document}

%% file: introduction.tex
Many scientific questions about networks amount to problems of {\em network
  comparison}: one wants to know whether networks observed at different times,
or in different locations, or under different environmental or experimental
conditions, actually differ in their structure.  Such problems arise in
neuroscience (e.g., comparing subjects with disease conditions to healthy
controls, or the same subject before and after learning), in biology (e.g.,
comparing gene- or protein- interaction networks across species, developmental
stages or cell types), and in social science (e.g., comparing different social
relations within the same group, or comparing social groups which differ in
some outcome).  That the graphs being compared are not identical or even
isomorphic is usually true, but scientifically unhelpful.  What we need is a
way to say if the difference between the graphs exceeds what we should expect
from mere population variability or stochastic fluctuations.  Network
comparison, then, is a kind of two-sample testing, where we want to know
whether the two samples could have come from the same source distribution.  It
is made challenging by the fact that the samples being compared are very
structured, high-dimensional objects (networks), and more challenging because
we often have only {\em one} graph in each sample.

We introduce a methodology for network comparison.  The crucial idea is to
approximate networks by continuous geometric objects, namely probability
densities, and then do two-sample bootstrap tests on those densities.
Specifically, we draw on recent work showing how many real-world networks are
naturally embedded in hyperbolic (negatively curved) manifolds.  Graphs then
correspond to clouds of points in hyperbolic space, and can be viewed as being
generated by sampling from an underlying density on that space.  We estimate a
separate density for each of the two networks being compared, calculate the
distance between those densities, and compare it to the distance expected under
sampling from a pooled density estimate.

Our method, while conceptually fairly straightforward, is admittedly more
complicated than the current practice in the scientific literature, which is to
compare networks by taking differences in {\em ad hoc} descriptive statistics
(e.g., average shortest path lengths, or degree distributions).  It is very
hard to assess the statistical significance of these differences, and
counter-examples are known where the usual summary statistics fail to
distinguish graphs which are qualitatively radically different (e.g., grid-like
graphs from highly clustered tree-like ones).  Similarly, whole-graph metrics
and similarity measures are of little statistical use, without probability
models to gauge their fluctuations.  Below, we show through simulations that
our method let us do network comparisons where (i) we can assess significance,
(ii) power is high for qualitative differences, and (iii) when we detect
differences, we also have some idea how {\em how} the networks differ.

%% file: motivation.tex
A fundamental issue with network comparison, mentioned in the introduction, is
that we often have {\em only} two networks to compare, and nonetheless need to
make some assessment of statistical significance.  This can obviously only be
done by regarding the networks as being drawn from (one or more) probability
models, and restricting the form of the model so that an observation of a
single graph is informative about the underlying distribution.  That is, we
must restrict ourselves to network models which obey some sort of law of large
numbers or ergodic theorem within a single graph, or else we always have $n=1$.
As in any other testing problem, the better the alignment between the model's
restrictions and actual properties of the graphs, the more efficiently the test
will use the available information.

\paragraph{Salient properties of actual networks} 
Over the last two decades, it has become clear that many networks encountered
in the real world, whether natural or human-made, possess a number of
mathematically striking properties \citep{MEJN-on-networks}.  They have highly
right-skewed degree distributions, they show the ``small-world effect'' of
short average path lengths (growing only logarithmically with the number of
nodes) but non-trivial transitivity of links, and high clusterability, often
with a hierarchical arrangement of clusters.  This is all a far cry from what
is expected of conventional random graphs.  While a large literature of
parametric stochastic models has developed to try to account for these
phenomena \citep{MEJN-on-networks}, there are few situations where a data
analyst can {\em confidently} assert that one of these models is even
approximately well-specified.

\paragraph{Current approaches to network comparison}

The typical approach in the literature is \emph{ad hoc} comparison of common
descriptive statistics on graphs (path lengths, clustering coefficients, etc.).
These statistics are often mis-applied, as in the numerous incorrect claims to
have found ``power law'' or ``scale-free'' networks \cite{power-law-or-not},
but that is not the fundamental issue.  Even the recent authoritative review
of, and advocacy for, the ``connectomics'' approach to neuroscience by Sporns
\cite{Sporns-networks-of-the-brain} takes this approach.  Disturbingly,
\cite{Henderson-Robinson-network-structure-in-cortex} shows that, with commonly
used choices of statistics and criteria, this approach cannot distinguish
between complex, hierarchically-structured networks, and simple two-dimensional
grids (such as a grid over the surface of the cortex).

More formally, \cite{Pao-Coppersmith-Priebe-power-of-graph-detection} study the
power of tests based on such summaries to detect departures from the null
hypothesis of completely independent and homogeneous edges (Erdos-Renyi graphs)
in the direction of independent but heterogeneous edges. Their results were
inconclusive, and neither the null nor the alternative models are plausible for
real-world networks.  Apart from this, essentially nothing is known about
either the significance of such comparisons or their power, how to combine
comparisons of different descriptive statistics, which statistics to use, or if
significant differences are found, how to infer changes in structure from them.
The issue of statistical significance also afflicts graph metrics and
similarity measures, even those with plausible rationales in graph theory
\citep{Koutra-Vogelstein-Faloutsos-deltacon}.

\cite{Hunter-Goodreau-Handcock-gof-of-social-networks} show one way to check
goodness-of-fit for a model of a single network, using simulations to check
whether the observed values of various graph statistics are plausible under the
model's sampling distribution.  But they are unable to combine checks with
different statistics, cannot find the power of such tests, and do not touch on
differences across networks.

More relevantly to comparisons, \cite{Inferring-network-mechanisms} use
machine-learning techniques to classify networks as coming from one or another
of various generative models, taking features of the network (such as the
counts of small sub-graphs, or ``motifs'') as the inputs to the classifier.
They demonstrate good operating characteristics in simulations, but rely on
having a good set of generative models to start with.

The approach to network comparison most similar to ours is
\cite{Rosvall-Bergstrom-mapping-change}, which like our proposed methods, uses
bootstrap resampling from models fit to the original networks to assess
significance of changes.  The goal there however is not to detect global
changes in the network structure, but local changes in which nodes are most
closely tied to one another.

\paragraph{Hyperbolic geometry of networks}
While waiting for scientifically-grounded parametric models, we seek a class of
non-parametric models which can accommodate the stylized facts of complex
networks.  Here we draw on the more recent observation that for many real-world
networks, if we view them as metric spaces with distance given by shortest path
lengths, the resulting geometry is {\em hyperbolic}
\citep{Albert-et-al-negative-curvature-of-networks,
  Kennedy-et-al-hyperbolicity-of-networks, Krioukov-et-al-hyperbolic-geometry},
rather than Euclidean.  Said another way, many real-world networks can be
naturally embedded into negatively-curved continuous spaces.  Indeed,
\citep{Krioukov-et-al-hyperbolic-geometry} show that if one draws points
representing nodes according to a ``quasi-uniform'' distribution on the
hyperbolic plane (see \eqref{eqn:parametrized.hyperbolic.densities} below), and
then connects nodes with a probability that decays according to the hyperbolic
distance between the points representing them, one naturally obtains graphs
showing right-skewed degree distributions, short average path lengths, and
high, hierarchical clusterability.

\paragraph{Continuous latent space models}
The model of \citep{Krioukov-et-al-hyperbolic-geometry} is an example of a {\em
  continuous latent space} model, characterized by a metric space
$(\ContinuousSpace,\ContinuousMetric)$, a link probability function $W$, and a
probability density $f$ on $\ContinuousSpace$, the {\em node density}.  Points
representing nodes are drawn iidly from $f$, and edges form
independently between nodes at $x$ and $y$ with probability $W(x,y) =
W(\ContinuousMetric(x,y))$ decreasing in the distance.  As a hierarchical model,
\begin{eqnarray}
\label{eqn:CLSM-as-hierarchical-model}
Z_i & \sim_{iid} & f\\
\nonumber A_{ij}| Z_1, \ldots Z_n & \sim_{ind} & W(\rho(Z_i,Z_j))
\end{eqnarray}
where $A_{ij}$ is the indicator variable for an edge between nodes $i$ and $j$.
Holding $\ContinuousSpace,\ContinuousMetric,W$ fixed, but allowing $f$ to vary,
we obtain different distributions over graphs.  Two densities $f,g$ on
$\ContinuousSpace$ determine the same distribution over graphs if $f$ is the
image of $g$ under some isometry of $(\ContinuousSpace,\ContinuousMetric)$.
Note that node {\em densities} can be compared regardless of the number of
nodes in the observed graphs.

The best-known continuous latent space model for social networks is that of
\citep{Hoff-Raftery-Handcock}, where the metric space is taken to be Euclidean
and the density $f$ is assumed Gaussian.  Our general methodology for network
comparison could certainly be used with such models.  However, the striking
properties of large real-world graphs, such as their highly-skewed degree
distributions, lead us to favor the sort of hyperbolic model used by
\citep{Krioukov-et-al-hyperbolic-geometry}, but without their restrictive
assumptions on $f$.  Rather, we will show how to non-parametrically estimate
the node density from a single observed graph, and then reduce network
comparison to a comparison of these probability densities.

Continuous latent space models are themselves special cases of models called
{\em graphons}, lifting the restriction that $\ContinuousSpace$ be a metric
space, and requiring of the edge probability function $W(x,y)$ only that it be
measurable and symmetric in its arguments\footnote{Graphons are often {\em
    defined} to have $\ContinuousSpace=[0,1]$ and $f$ Lebesgue measure.  One
  can show that any graphon over another measure space or with another node
  density is equivalent to one of this form, i.e., generates the same
  distribution over infinite graphs \citep[ch.\ 7]{Kallenberg-symmetries}.}.
Any distribution over infinite graphs which is invariant under permuting the
order of the nodes turns out to be a mixture of such graphons \citep[ch.\
7]{Kallenberg-symmetries}.  Moreover, as one considers larger and larger
graphs, the properties of the observed graph uniquely identify the generating
graphon \citep{Diaconis-Janson-graph-limits}; what almost comes to the same
thing, the limit of a sequence of growing graphs {\em is} a graphon
\citep{Borgs-Chayes-Lovasz-et-al-graph-limits-and-parameter-testing,
  Lovasz-large-networks, Borgs-Chayes-Zhao-sparse-graph-convergence}.  One
might, then, try to use our approach to compare graphons with estimated $f$ and
$W$.  While graphon estimation is known to be possible in principle
\citep{Bickel-Chen-Levina-method-of-moments,
  Choi-Wolfe-consistency-of-co-clustering}, there are no published,
computationally feasible methods to do it.  Moreover, we expect to gain power
by tailoring our models to enforce salient network properties, as described
above.  Accordingly, we turn to some of the important aspects of hyperbolic
geometry.

\subsection{Hyperbolic spaces}
\label{subsec:hyperbolic spaces}

Hyperbolic spaces are metric spaces which are negatively curved --- the angles
in a triangle of geodesics sum to less than $180$ degrees.  The oldest example
of such a space is the surface of the hyperboloid, the surface of points $(x_1,
x_2,x_3) \in \R^3$ such that 
\begin{equation*}
x_{1}^2+x_{2}^2-x_{3}^2=1,
\end{equation*}
with the distance between points taken to be the Euclidean length of the
shortest path between them along the surface.  Another, and perhaps even more
basic, example of a hyperbolic space is a tree, again with the shortest-path
metric.  Our starting data will be observed networks, which are typically at
least locally tree-like, and so also possess a hyperbolic geometry
\citep{Jonckheere-et-al-scaled-hyperbolic-graphs}.

As explained above, we aim to represent this discrete hyperbolic geometry with
a density over a continuous hyperbolic space.  For concreteness, we will focus
on the hyperbolic plane $\hyperboloid{2}$, whose most basic geometric model is
just the surface of the hyperboloid.  It will be more convenient to work with
another model of $\hyperboloid{2}$: the {\em Poincar{\'e} half-plane} of $\C$,
$$\hyperboloid{2}=\left\{x+iy \mid x\in\mathbb{R}, y\in(0,\infty)\right\}$$
equipped with the metric $d\ContinuousMetric^2=(dx^2+dy^2)/y^{2}$.

\begin{figure}
  \begin{center}
    \includegraphics[width=40mm,height=40mm]{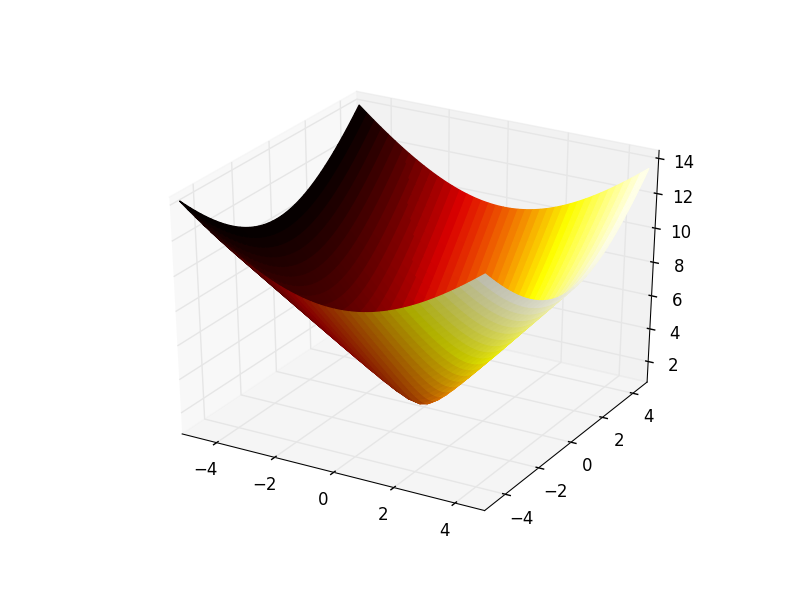}
    \includegraphics[width=30mm,height=20mm]{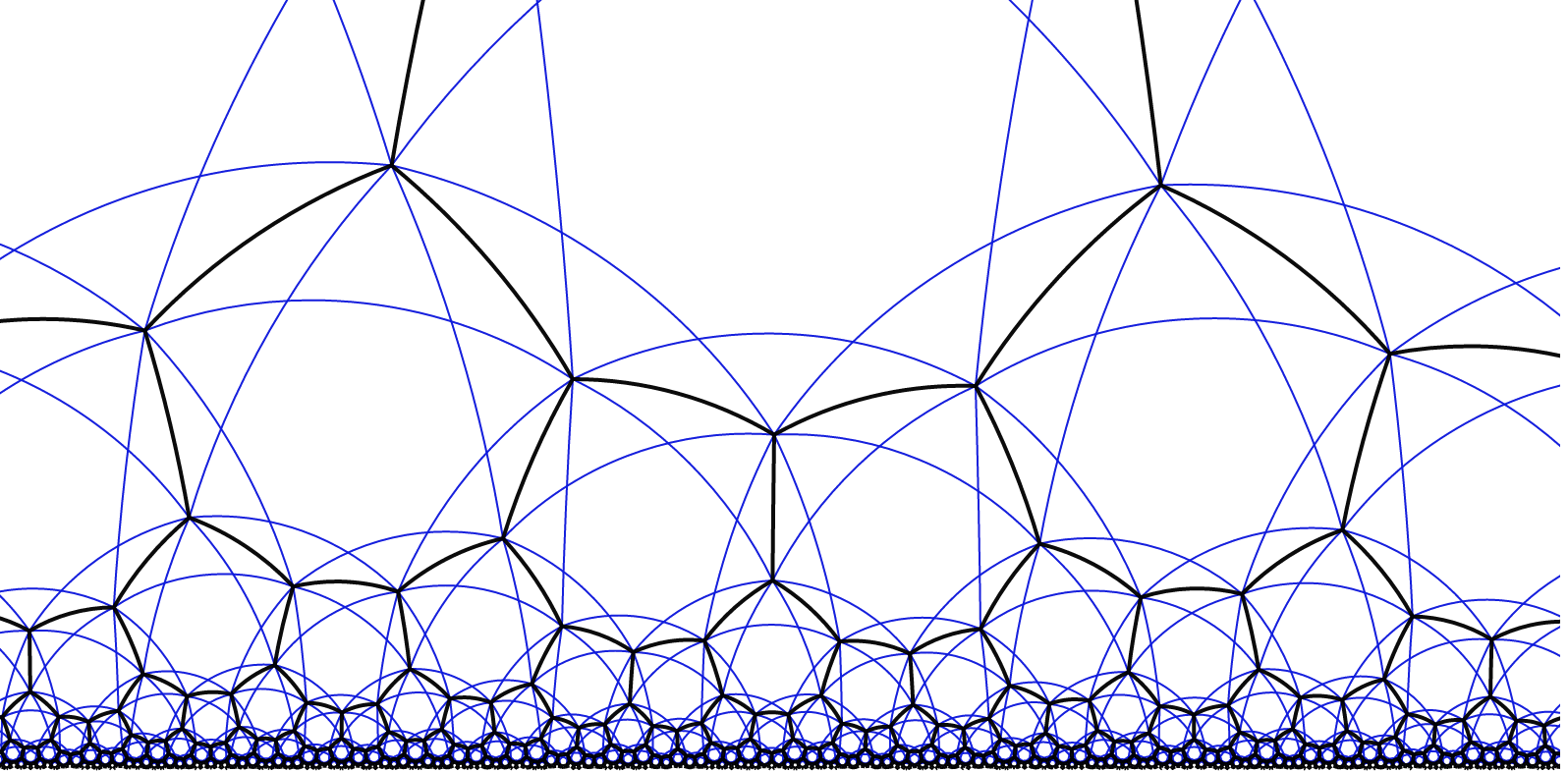}
  \end{center}
  \caption{{\bf Models of $\hyperboloid{2}$} A connected component of the
    hyperboloid $x_{3}^2=1+x_{1}^2+x_{2}^2$ (left), with the metric given by
    the Euclidean length of the shortest path between points along the surface,
    is isometric to the {\em Poincar{\'e} half-plane} (right) under a suitable
    non-Euclidean metric.  The half-plane is tiled into regions of equal area
    with respect to the metric. (Images from \citep{halfplane}, under a
    Creative Commons license.)}
  \label{fig:H2}
\end{figure}

As mentioned above, \citep{Krioukov-et-al-hyperbolic-geometry} showed that if
the density of nodes on the Poincar{\'e} half-plane is one of the {\em
  quasi-uniform} densities,
\begin{equation}
  \label{eqn:parametrized.hyperbolic.densities}
  \QuasiUniformDensity_{\delta,R}(re^{i\theta})=\frac{\delta\sinh{\delta r}}{2\pi(\sinh{r})\cosh{(\delta R-1)}},\quad\delta>0
\end{equation}
one obtains graphs which reproduce the stylized facts of right-skewed degree
distributions, clusterability, etc., for a wide range of link probability
functions $W$, including Heaviside step functions
$\Theta(\ContinuousMetric-c)$.  Note that the mode of $\QuasiUniformDensity$ is
always at $0+i$, with the parameter $\delta > 0$ controlling the dispersion
around the mode, and $R >0$ being an over-all scale factor.  As $\delta$ grows,
the resulting graphs become more clustered.

We will introduce another family of densities on $\hyperboloid{2}$, the
hyperboloid \textit{Gaussians}, in the next section.

\begin{figure}
  \begin{center}
    \includegraphics[width=25mm,height=25mm]{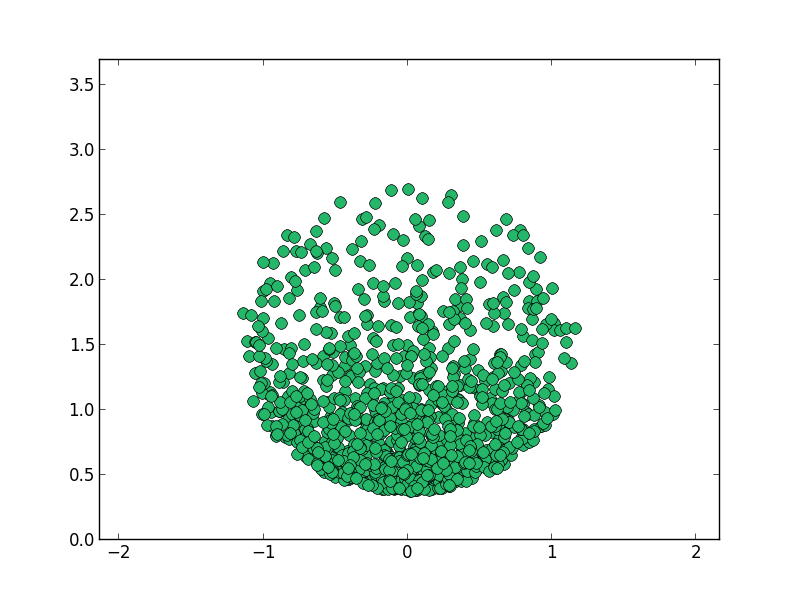}
    \includegraphics[width=25mm,height=25mm]{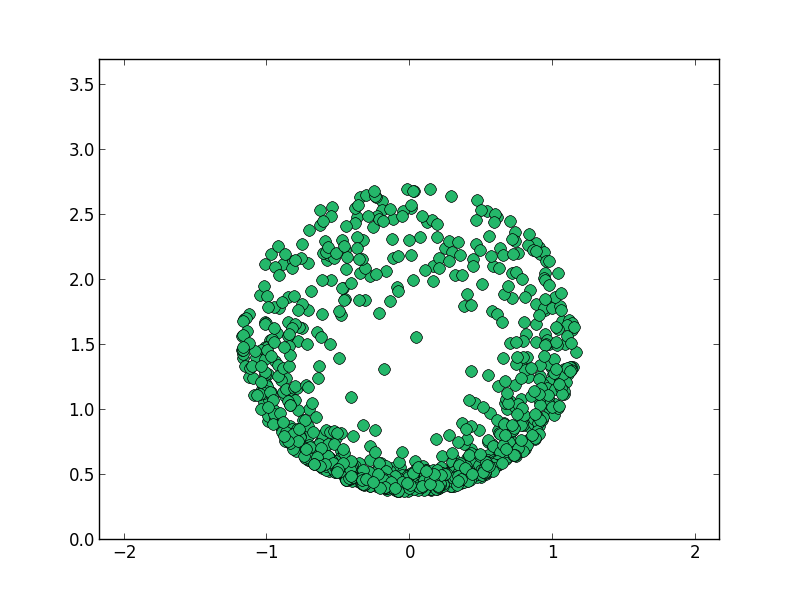}
    \includegraphics[width=25mm,height=25mm]{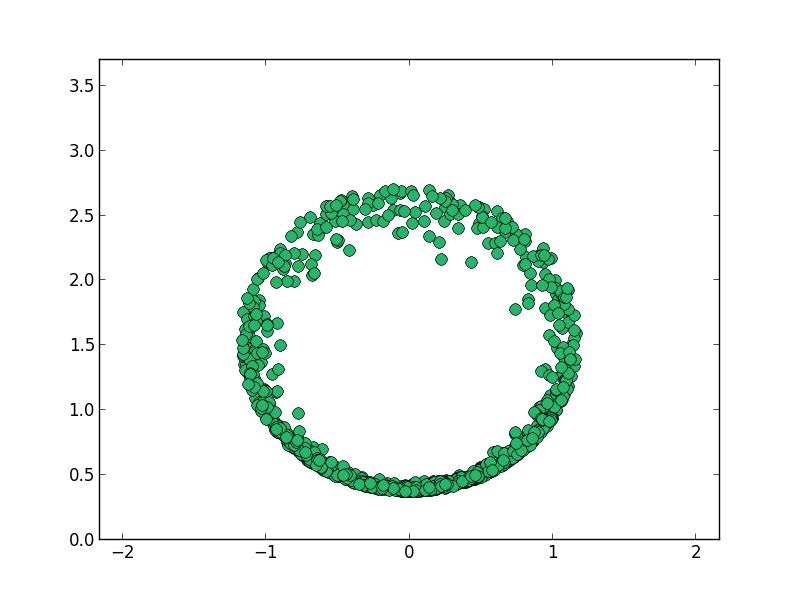}
    \includegraphics[width=25mm,height=25mm]{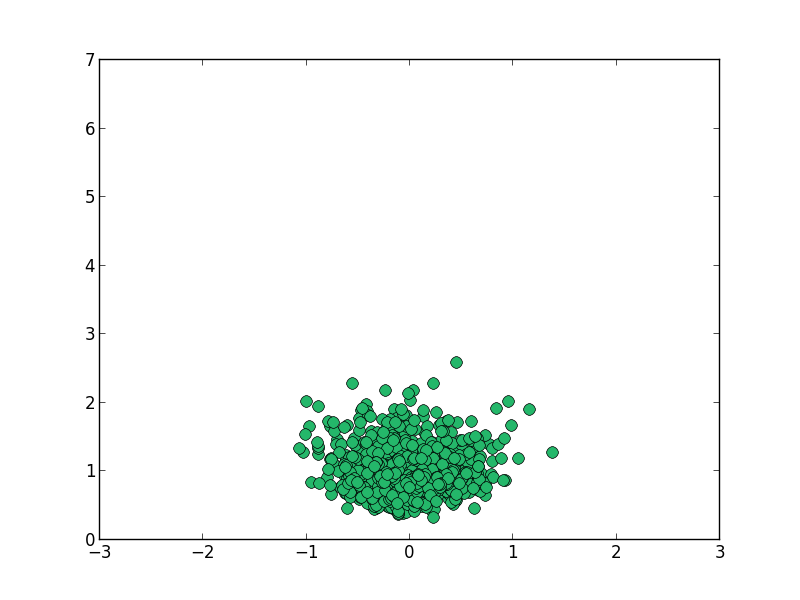}
    \includegraphics[width=25mm,height=25mm]{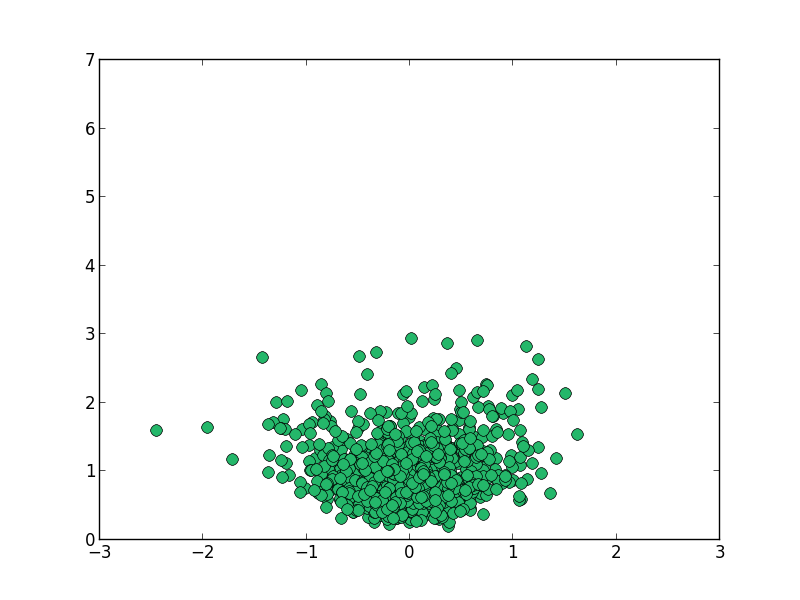}
    \includegraphics[width=25mm,height=25mm]{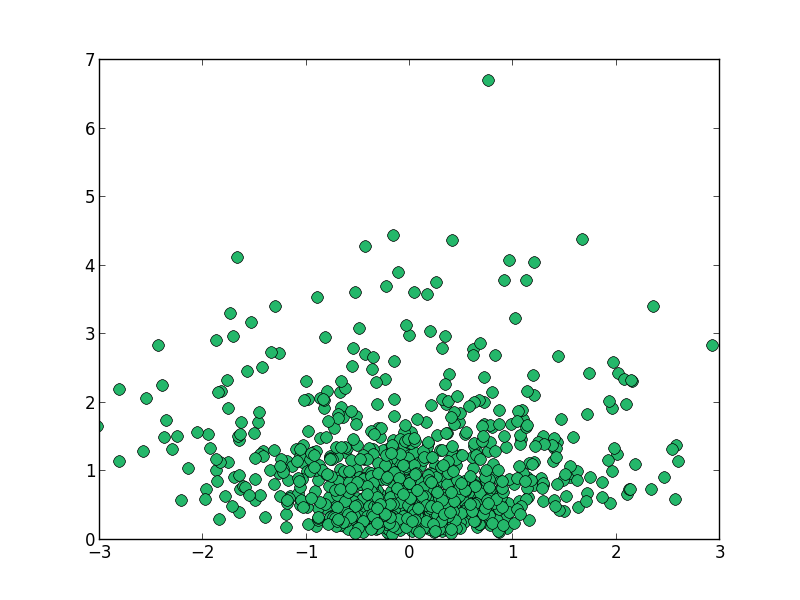}
  \end{center}
  \caption{{\bf Densities on $\hyperboloid{2}$} 1000 points drawn iidly from
    quasi-uniform densities, Eq.\ \ref{eqn:parametrized.hyperbolic.densities}
    (top; $\delta=1,10,30$ from left to right, $R=1$ throughout), and from
    hyperbolic Gaussian densities, Eq.\ \ref{eqn:hyperbolic-gaussian-density}
    (bottom, $\HyperbolicGaussianParam=0.05,0.1,0.3$ from left to right).}
  \label{fig:quasi-uniform.densities}
\end{figure}

\begin{figure}
  \begin{center}
    \includegraphics[width=25mm,height=25mm]{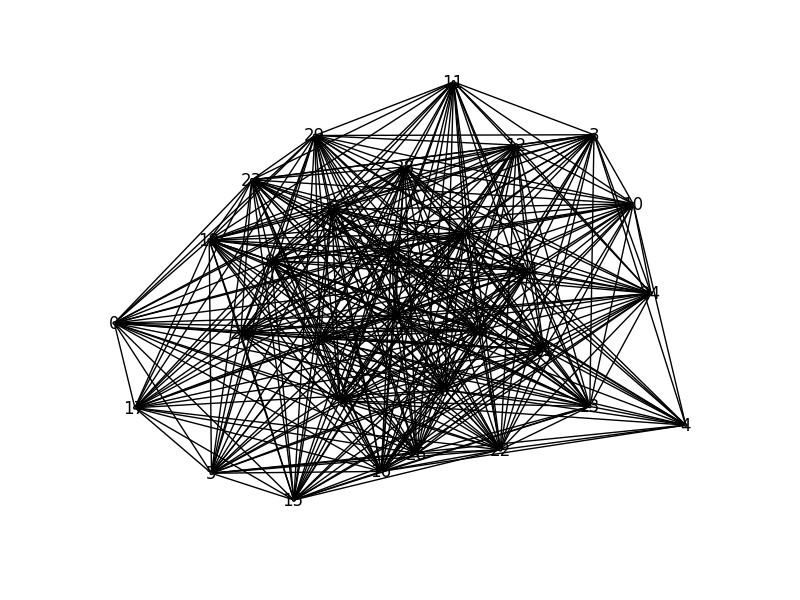}
    \includegraphics[width=25mm,height=25mm]{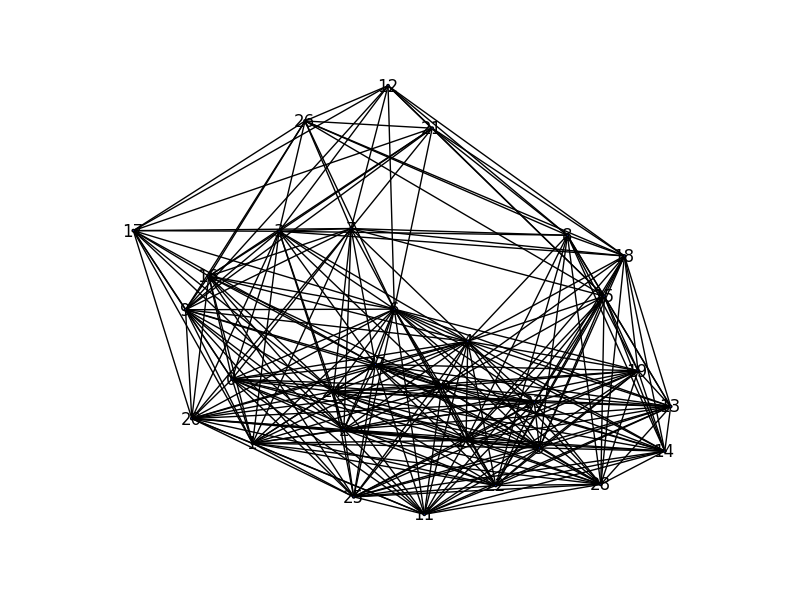}
    \includegraphics[width=25mm,height=25mm]{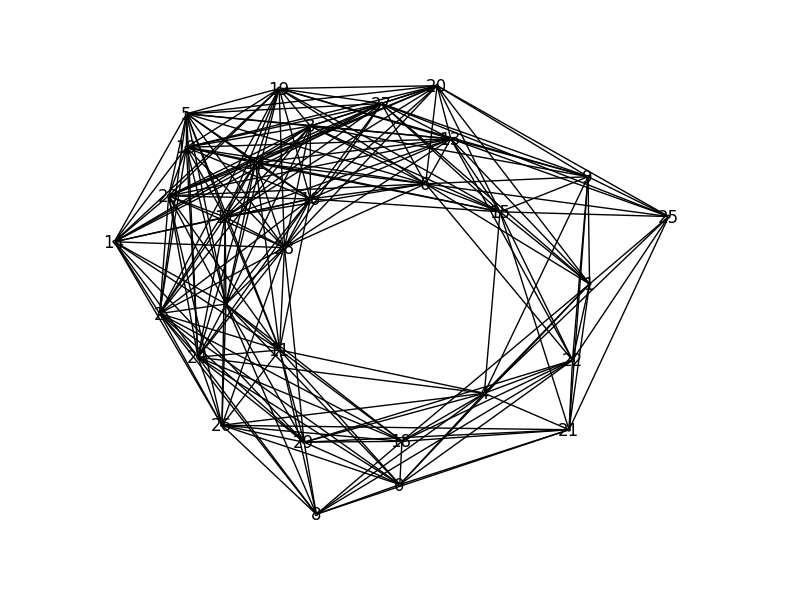}
    \includegraphics[width=25mm,height=25mm]{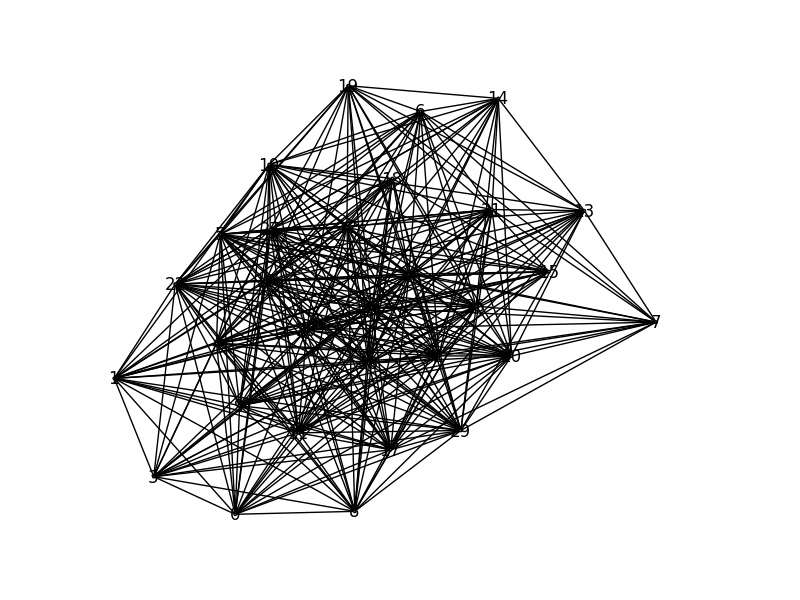}
    \includegraphics[width=25mm,height=25mm]{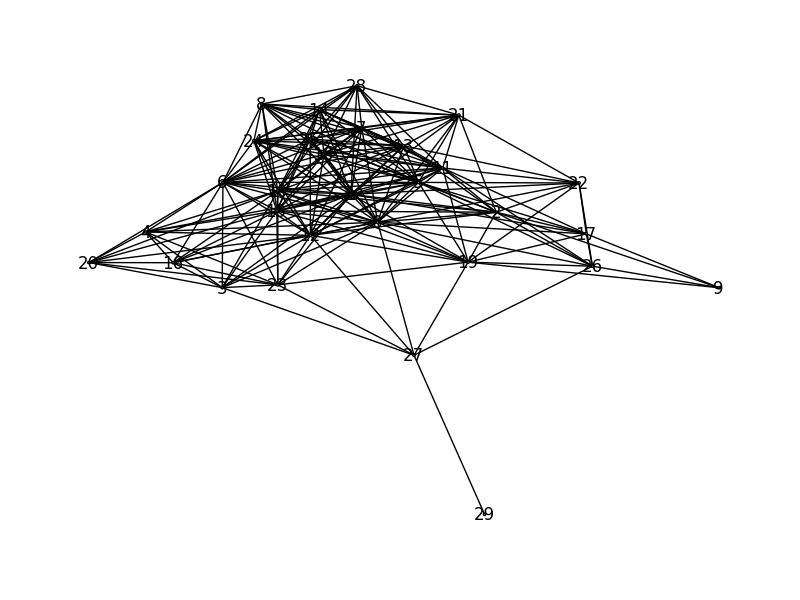}
    \includegraphics[width=25mm,height=25mm]{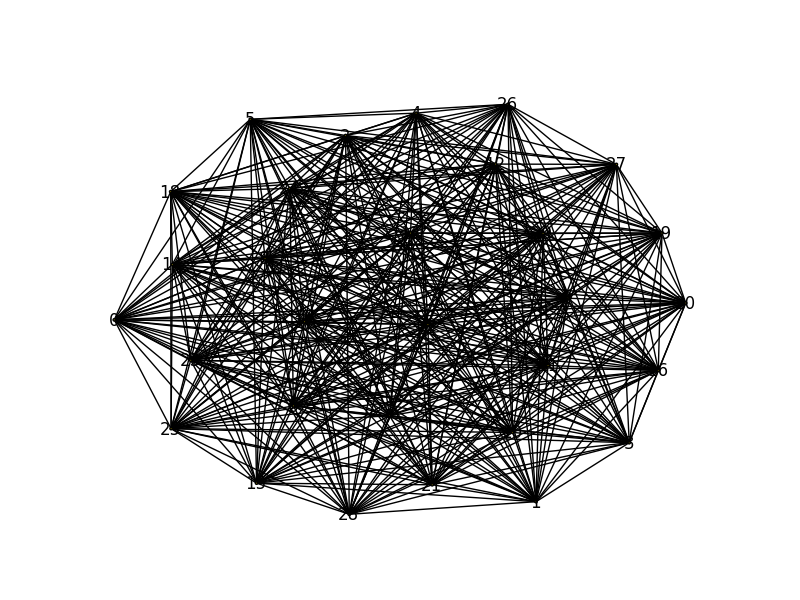}
  \end{center}
  \caption{{\bf Hyperbolic latent-space graphs} Graphs formed by drawing 30
    node locations as in Fig.\ \ref{fig:quasi-uniform.densities}, and applying
    the link probability function $W(x,y) = \Theta(\ContinuousMetric(x,y) -
    1.5)$.  Note how the graphs in the bottom row become more clustered as the
    $\delta$ parameter increases from left to right.}
  \label{fig:generated.graphs}
\end{figure}

Fig.\ \ref{fig:quasi-uniform.densities} shows samples from quasi-uniform
distributions on $\hyperboloid{2}$, and Fig.\ \ref{fig:generated.graphs} the
resulting graphs.  While we will use such networks as test cases, we emphasize
that we will go beyond \eqref{eqn:parametrized.hyperbolic.densities} to a fully
nonparametric estimation of the node density.

%% file: methodology.tex
\begin{figure}
  \begin{center}
    \includegraphics[width=70mm,height=50mm]{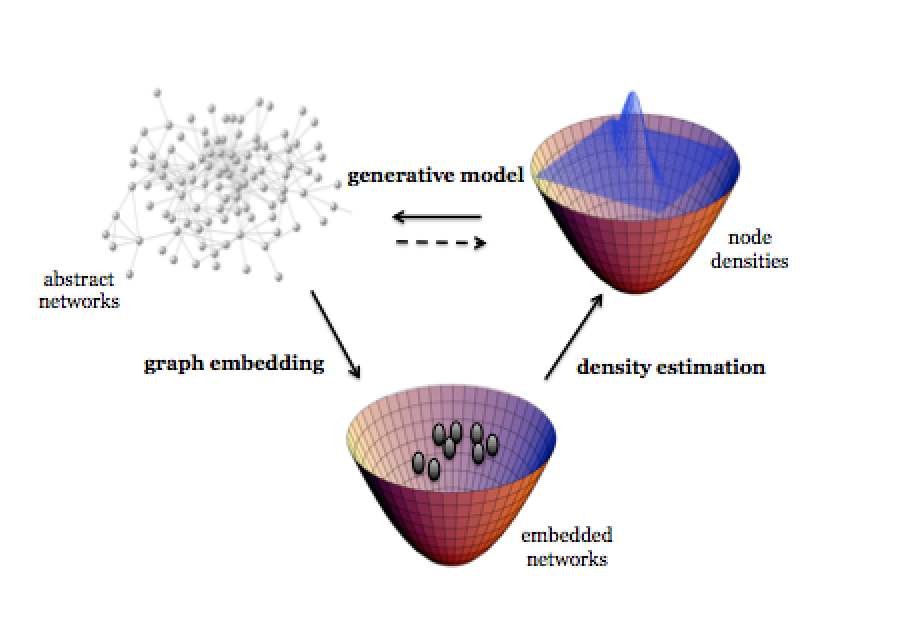}
  \end{center}
  \caption{Schematic of network inference}
  \label{fig:methodology}
\end{figure}

Our goal is to compare networks by comparing node densities.  Our procedure for
estimating node densities has in turn two steps (Figure \ref{fig:methodology}):
we embed the nodes of an observed network into $\hyperboloid{2}$ (\S
\ref{sec:embedding}), and then estimate a density from the embedded
points (\S \ref{sec:gkde}).  We may then compare the observed difference
between estimated node densities from two graphs to what would be expected if
we observed two graphs drawn from a common node density (\S
\ref{sec:comparison}).

\subsection{Graph embedding}
\label{sec:embedding}

An {\em embedding} of a graph $G$ is a mapping of its nodes $V_G$ to points
into a continuous metric space $(\ContinuousSpace,\ContinuousMetric)$ which
preserves the structure of the graph, or tries to.  Specifically, the distances
between the representative points should match the shortest-path distances
between the nodes, as nearly as possible.  This is a multidimensional scaling
problem, where typically one seeks the embedding $\phi: V_G \mapsto
\ContinuousSpace$ minimizing
\begin{equation}
  \label{eqn:mds}
  \sum_{(v,w)\in V_G^2}(\rho_G(v,w)-\ContinuousMetric(\phi(v),\phi(w)))^2,
\end{equation}
where $\rho_G$ is the shortest-path-length metric on $V_G$.  Classically, when
$\ContinuousSpace=\R^n$ and $\ContinuousMetric$ is the Euclidean metric, the
arg-min of \eqref{eqn:mds} can be found by spectral decomposition of the matrix
of $\rho_G(v,w)$ values \citep[ch.\ 3]{Hand-Mannila-Smyth}.

Spectral decomposition does not however give the arg-min of \eqref{eqn:mds}
when $\ContinuousSpace=\hyperboloid{2}$ with the appropriate non-Euclidean
metric.  While the solution could be approximated by gradient descent
\citep{Cvetkovski-Crovella-mds-on-the-Poincare-disk}, we follow
\citep{Begelfor-Werman-world-is-not-always-flat} in changing the problem
slightly.  They propose minimizing
\begin{equation}
  \label{eqn:general.mds}
  \sum_{(v,w)\in V_G^2}{(\cosh{\rho_G(v,w)}-\cosh{\ContinuousMetric(\phi(v),\phi(w))})^2}
\end{equation}
which can be done exactly via a spectral decomposition.  Specifically, let
$R_{ij} = \cosh{\rho_G(i,j)}$, whose leading eigenvector is $u_1$ and whose
trailing eigenvectors are $u_2$ and $u_3$.  Then the $i^{\mathrm{th}}$ row of
the matrix $(u_1 u_2 u_3)$ gives the $\hyperboloid{2}$ coordinates for node
$i$.  If $R$ has one positive eigenvalue, exactly 2 negative eigenvalues, and
all remaining eigenvalues vanish, this defines an exact isometric embedding
\citep{Begelfor-Werman-world-is-not-always-flat}.

We have not found a way of estimating the node density which avoids the initial
step of embedding.  Our method is, however, fairly indifferent as to {\em how}
the nodes are embedded, so long as this is done well, and in a way which does
not pre-judge the form of the node density.

\begin{figure}[h]
  \begin{center}
    \includegraphics[width=25mm,height=25mm]{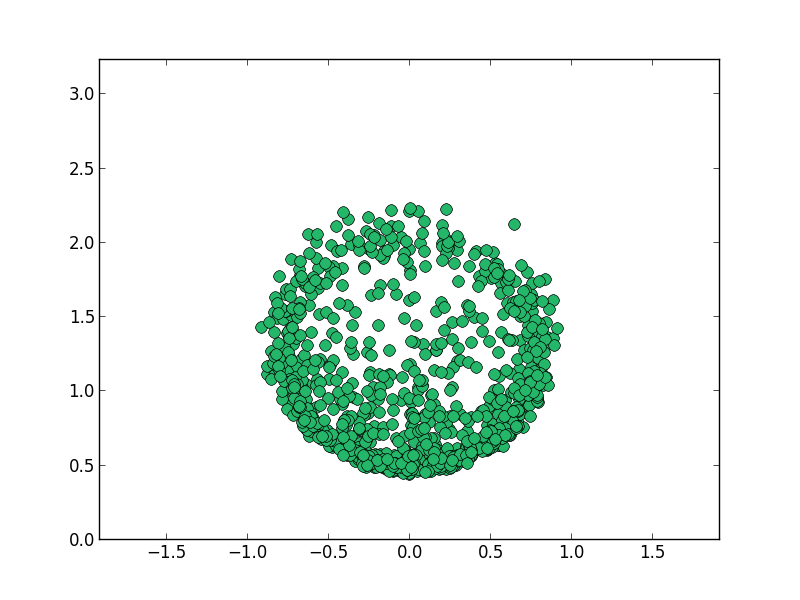}
    \includegraphics[width=25mm,height=25mm]{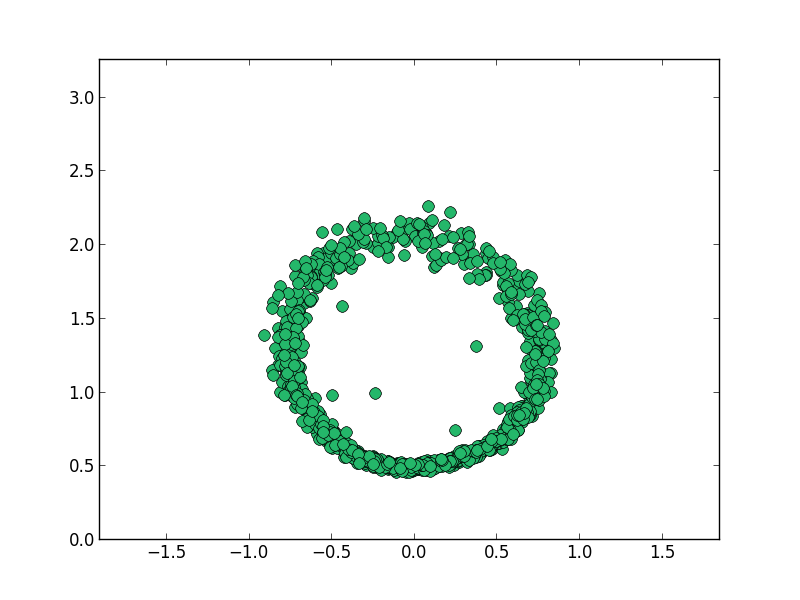}
    \includegraphics[width=25mm,height=25mm]{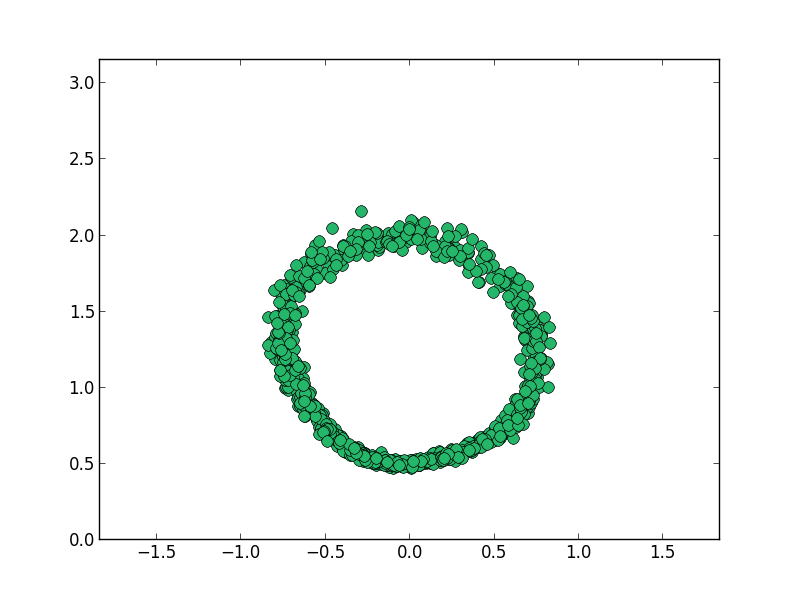}
    \includegraphics[width=25mm,height=25mm]{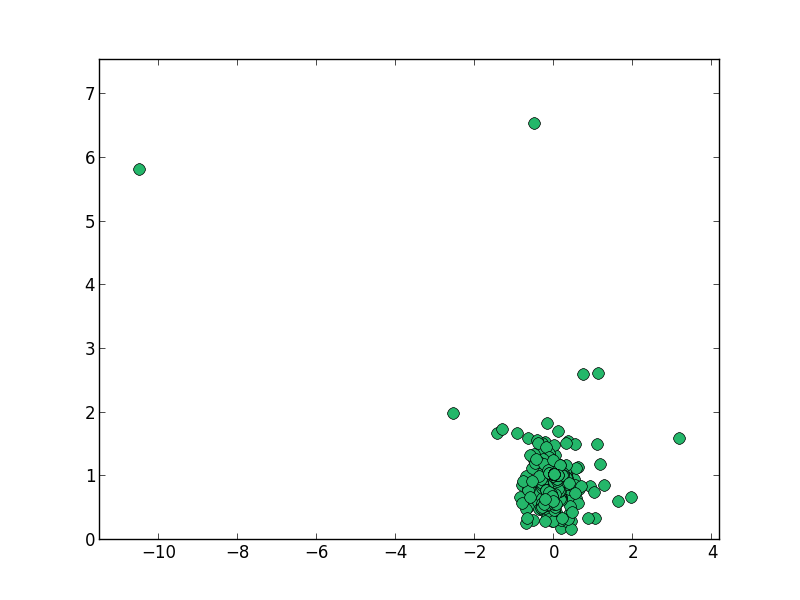}
    \includegraphics[width=25mm,height=25mm]{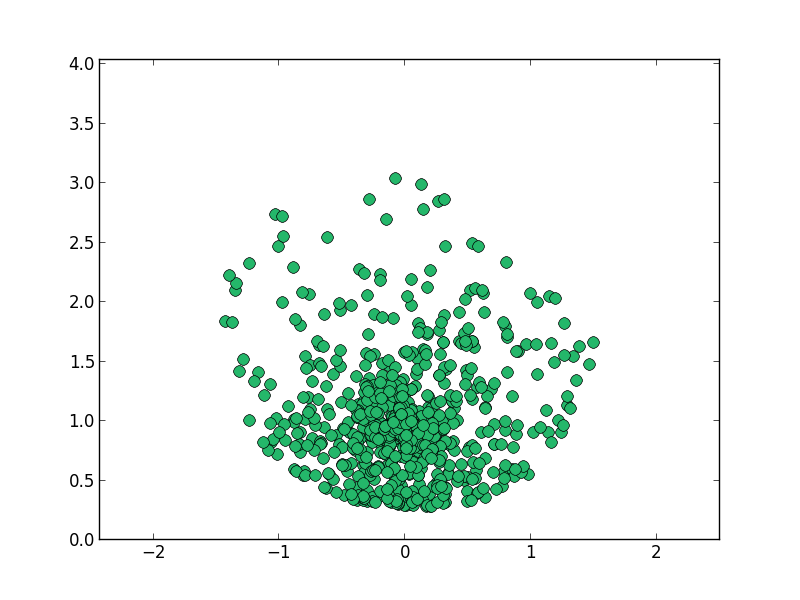}
    \includegraphics[width=25mm,height=25mm]{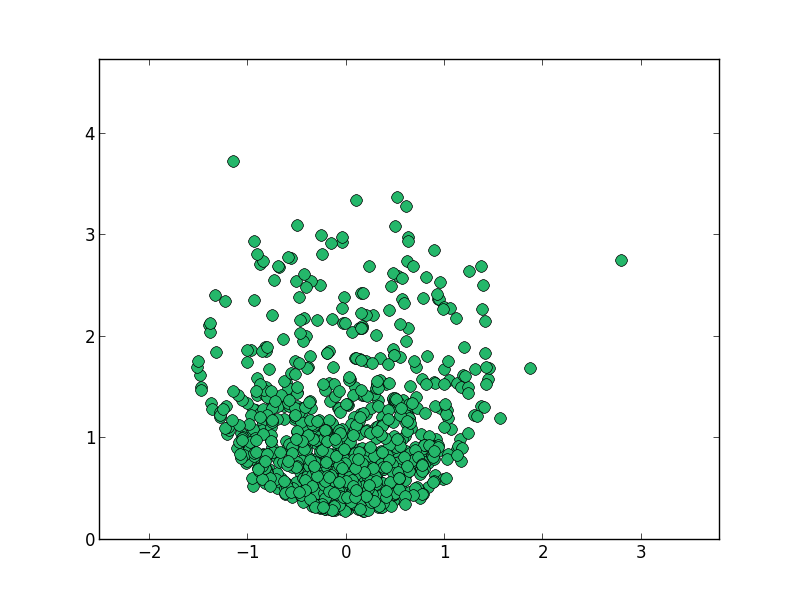}
  \end{center}
  \caption{{\bf Reembedded Generated Graphs} Results of embedding simulated
    graphs, formed as in Fig. \ref{fig:generated.graphs}, back into
    $\hyperboloid{2}$.  Comparison with Fig. \ref{fig:quasi-uniform.densities}
    illustrates the fidelity of the embedding process.}
  \label{fig:mds}
\end{figure}

\subsection{Density estimation}
\label{sec:gkde}

Having embedded the graph into $\hyperboloid{2}$, we estimate the node density.
Our procedure for doing so is more easily grasped by first reviewing the
connections between kernel density estimation, convolution, and Fourier
transforms in Euclidean space.

\paragraph{Kernel Density Estimation in Euclidean Space as Convolution}
In Euclidean space, kernel density estimation smooths out the empirical
distribution by adding a little bit of noise around each observation.  Given
observations $z_1, z_2, \ldots z_n \in \R^p$, and a normalized kernel function
$K_h$, the ordinary kernel density estimator $\KDE$ at a point $z \in \R^p$ is
\begin{eqnarray*}
\KDE(z) & = & \frac{1}{n}\sum_{i=1}^{n}{K_h(z-z_i)}\\
& = & \int_{\R^p}{K_h(z-z^{\prime}) \left(\frac{1}{n}\sum_{i=1}^{n}{\delta(z^{\prime}-z_i)}\right) dz^{\prime}}\\
& = & \int_{\R^p}{K_h(z-z^{\prime}) \EmpiricalMeasure(dz^{\prime})}\\
& = & (K_h * \EmpiricalMeasure)(z)
\end{eqnarray*}
where the third line defines the empirical measure $\EmpiricalMeasure$, and $*$
denotes convolution.  In words, the kernel density estimate is the convolution
of the empirical measure with the kernel.  Here the role of the kernel $K_h$ is
not so much to be a distribution over the Euclidean space, as a distribution
over {\em translations} of the space: $K_h(z-z_i)$ is really the density at the
translation mapping the data point $z_i$ into the operating point $z$.  As it
happens, the group of translations of $\R^p$ is also $\R^p$, but when we adapt
to non-Euclidean spaces, this simplifying coincidence goes away.

Since, in Euclidean space, the Fourier transform $\mathcal{F}$ converts
convolutions into products \citep{stein1971introduction},
\[
\Fourier{\KDE}(s) = \Fourier{K_h}(s) \Fourier{\EmpiricalMeasure}(s)
\]
This relation often greatly simplifies computing $\KDE$.  It also lets us
define the bandwidth $h$, through the relation
$\Fourier{K_h}(s)=\Fourier{K}(hs)$. 

It is well known that kernel density estimators on $\R^p$, with $h \rightarrow
0$ at the appropriate rate in $n$, are minimax-optimal in their $L_2$ risk
\citep{van-der-Vaart-asymptotic-stats}.  With suitable modifications, this
still holds for compact manifolds \citep{Pelletier-KDE-on-manifolds}, but the
hyperbolic plane $\hyperboloid{2}$ is not compact.

\subsubsection{$\hyperboloid{2}$-Kernel density estimator}

Our method for density estimation on $\hyperboloid{2}$ is a generalization of
Euclidean kernel density estimation.  In $\R^p$, the kernel is a density on
translations of $\R^p$.  For $\hyperboloid{2}$, the appropriate set of
isometric transformations are not translations, but rather the class of
``M{\"o}bius transformations'' represented by the Lie group $\SL_2$
\citep{Terras-harmonic-1, Huckemann-et-al-mobius-deconvolution}.  An
$\hyperboloid{2}$ kernel, then, is a probability density on $\SL_2$.  We may
write $K_h(z,z_i)$ to abbreviate the density the kernel $K_h$ assigns to the
M{\"o}bius transform taking $z_i$ to $z$.  
The generalized kernel density estimator on $\hyperboloid{2}$ takes the form
\begin{eqnarray}
  \label{eqn:H2.estimator}
  \gKDE(z) & = & \frac{1}{n}\sum_{i=1}^n{K_h(z,z_i)}\\
  & = & (K_h * \EmpiricalMeasure)(z)
\end{eqnarray}

In Euclidean space, the Fourier transform analyzes functions (or generalized
functions, like $\widehat{P}_n$) into linear combinations of the eigenfunctions
of the Laplacian operator.  The corresponding operation for $\hyperboloid{2}$
is the {\em Helgason}, or {\em Helgason-Fourier}, transform $\mathcal{H}$
\citep{Terras-harmonic-1}.  The Fourier basis functions are indexed by $\R^p$,
which is the group of translations; for analogous reasons, the Helgason basis
functions are indexed by $\C \times \SO_2$.  Many of the formal properties of
the Fourier transform carry over to the Helgason transform.  (See App.\
\ref{app:helgason}.)  In particular, convolution still turns into
multiplication:
\begin{equation}
\label{eqn:helgason-transform-of-gkde}
\Helgason{\gKDE} = \Helgason{K_h} \Helgason{\EmpiricalMeasure},
\end{equation}
where $\mathcal{H}[K_h]$ denotes the Helgason-Fourier transform of the well-defined density on $\mathbb{H}_2$ induced by the density $K_h$ on $\SL_2$, and we define the bandwidth $h$ through
$$\Helgason{K_h}(s,M)=\Helgason{K}(hs,M).$$
As in Euclidean density estimation, $h$ may be set through cross-validation.

In a separate manuscript \citep{Asta-gKDE}, we show that the $L_2$ risk of
\eqref{eqn:H2.estimator} goes to zero at the minimax-optimal rate, under mild
assumptions on the smoothness of the true density, and of the kernel $K$.
(This is a special case of broader results about generalized kernel density
estimation on symmetric spaces.) The assumptions on the kernel are satisfied by
what \citep{Huckemann-et-al-mobius-deconvolution} calls ``hyperbolic
Gaussians'', densities on $\hyperboloid{2}$ with parameter $\rho$ defined
through their Helgason transforms,
\begin{equation}
\label{eqn:hyperbolic-gaussian-density}
\Helgason{K}(s,M)\propto e^{\rho\overline{s(s-1)}} ~.
\end{equation}
Just as the ordinary Gaussian density is the unique solution to the heat
equation in Euclidean space, the hyperbolic Gaussian is the unique
($\mathbb{SO}_2$-invariant) solution to the heat equation on $\hyperboloid{2}$
\citep{Terras-harmonic-1}.  

\subsection{Network comparison}
\label{sec:comparison}

Combining embedding with kernel density estimation in $\hyperboloid{2}$ gives
us a method of estimating node densities, and so of estimating a hyperbolic
latent space model for a given network.  We now turn to {\em comparing}
networks, by comparing these estimated node densities.

Our method follows the general strategy advocated in
\citep{Network-comparisons}.  Given two graphs $G_1$ and $G_2$, we may
estimate two separate network models
$$\EstModel_1=\EstModel(G_1),\quad\widehat{P}_2=\widehat{P}(G_2).$$
We may also pool the data from the two
graphs to estimate a common model
$$\EstModel_{12} = \EstModel(G_1,G_2).$$
We calculate a distance $d^* = d(\EstModel_1,\EstModel_2)$ using any suitable
divergence.  We then compare $d^*$ to the distribution of distances which may
be expected under the pooled model $\EstModel_{12}$.  To do so, we
independently generate $G^{\prime}_1, G^{\prime}_2 \sim \EstModel_{12}$, and
calculate
$$d(\EstModel(G^{\prime}_1),\EstModel(G^{\prime}_2)).$$
That is, we bootstrap two independent graphs out of the pooled model, fit a
model to each bootstrapped graph, and calculate the distance between them.
Repeated over many bootstrap replicates, we obtain the sampling distribution of
$d$ under the null hypothesis that $G_1$ and $G_2$ are drawn from the same
source, and any differences between them are due to population variability or
stochastic fluctuations.\footnote{This method extends easily to comparing sets
  of graphs, $G_{11}, G_{12}, \ldots G_{1n}$ vs.\ $G_{21}, G_{22}, \ldots
  G{_{2m}}$, but the notation grows cumbersome.}

In our case, we have already explained how to find $\EstModel_1$ and
$\EstModel_2$.  Since we hold the latent space $\ContinuousSpace$ fixed at
$\hyperboloid{2}$, and the link probability function $W$ fixed, we can label
our models by their node densities, $\gKDE_1$ and $\gKDE_2$.  To obtain the
pooled model $\EstModel_{12}$, we first embed $G_1$ and $G_2$ separately using
generalized multidimensional scaling, and then do kernel density estimation on
the union of their embedded points.

The generalized multidimensional scaling technique we use depends only on the
eigendecomposition of matrices determined by shortest path lengths.
Therefore the $L_2$ difference
\begin{equation}
  \label{eqn:test}
  \|\gKDE_1-\gKDE_2\|_2
\end{equation}
between two estimated node densities $\gKDE_1,\gKDE_2$ is $0$ if and only
if the original sets of vertices from the different samples are isometric and hence (\ref{eqn:test}) approximates a well-defined metric $d$ on our continuous latent space models.  Moreover, since the
Plancherel identity carries over to the Helgason-Fourier transform
\citep{Terras-harmonic-1},
\begin{equation}
  \label{eqn:plancharel}
  \NodeDensityLpMetric{2}(f_1,f_2)=\|\Helgason{f_1}-\Helgason{f_2}\|_2,
\end{equation}
and, for our estimated node densities, $\Helgason{f}$ is given by
\eqref{eqn:helgason-transform-of-gkde}.  Appendix\
\ref{app:test-static-computing} gives full details on our procedure for
computing the test statistic \eqref{eqn:plancharel}.

\subsection{Theoretical Considerations}

Let us sum up our method, before turning to theoretical considerations.  (0) We
observe two graphs, $G_1$ and $G_2$.  (1) Through multi-dimensional scaling, we
embed them separately in $\hyperboloid{2}$ (\S \ref{sec:embedding}), getting
two point clouds, say $\mathbf{Z}_1$ and $\mathbf{Z}_2$.  (2) From each cloud,
we estimate a probability density on $\hyperboloid{2}$, using hyperbolic
Gaussian kernels, getting $\gKDEA$ and $\gKDEB$.  (\S \ref{sec:gkde}.)  We
calculate $\|\gKDEA-\gKDEB\|_2$ using \eqref{eqn:plancharel}.  We also form a
third density estimate, $\gKDEAB$, from $\mathbf{Z}_1 \cup \mathbf{Z}_2$.  (3)
We generate two independent graphs $G_1^*, G_2^*$ from $\gKDEAB$ according to
\eqref{eqn:CLSM-as-hierarchical-model}, and subject these graphs to
re-embedding and density estimation, obtaining $\gKDEBootA$ and $\gKDEBootB$
and so $\|\gKDEBootA - \gKDEBootB \|_2$.  Finally, (4) repeating step (3) many
times gives us the sampling distribution of the test statistic under the null
hypothesis that $G_1$ and $G_2$ came from the same source, and the $p$-value is
the quantile of $\|\gKDEA-\gKDEB\|_2$ in this distribution.

The final step of computing the $p$-value is a fairly unproblematic bootstrap
test.  The previous step of generating new graphs from the pooled model is also
an unproblematic example of a model-based bootstrap.  The kernel density
estimates themselves are consistent, and indeed converge at the minimax rate
\citep{Asta-gKDE}, {\em given} the point clouds on the hyperbolic plane.
This makes it seem that the key step is the initial embedding.  Certainly, it
would be convenient if the graphs $G_1$ and $G_2$ were generated by a
hyperbolic latent space model, and the embedding was a consistent estimator of
the latent node locations.  However, such strong conditions are not {\em
  necessary}.  Suppose that if $G_1 \sim \ProbModel_1$ and $G_2 \sim
\ProbModel_2 \neq \ProbModel_1$, then $\gKDE_1 \rightarrow f_1$ and $\gKDE_2
\rightarrow f_2$, with $\| f_1 - f_2\|_2 > 0$.  Then at any nominal size
$\alpha > 0$, the power of the test will go to 1.  For the nominal size of the
test to match the actual size, however, will presumably require a closer
alignment between the hyperbolic latent space model and the actual generating
distribution.

%% file: experiments.tex
\begin{figure}
  \begin{center}
    \includegraphics[width=70mm,height=50mm]{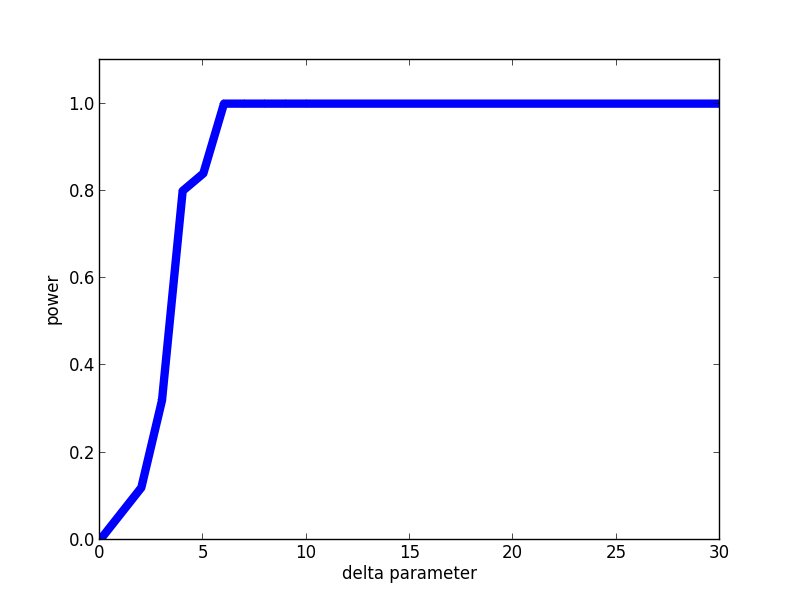}
  \end{center}
  \caption{{\bf Comparing Quasi-Uniforms} Power of our test, at size $\alpha =
    0.1$, for detecting the difference between a 100-node graph generated from
    the quasi-uniform density $q_{1,1}$ and a 100-node graph generated from
    $q_{\delta,1}$, as a function of the dispersion parameter $\delta$.}
  \label{fig:test}
\end{figure}

\paragraph{Comparison of Graphs with Quasi-Uniform Node Densities}
In our first set of simulation studies, we generated graphs which exactly
conformed to the hyperbolic latent space model, and in fact ones where the node
density was quasi-uniform (as in Fig.\ \ref{fig:generated.graphs}, top).  One
graph had 100 nodes, with latent locations drawn from a $q_{1,1}$ distribution;
the other, also of 100 nodes, followed a $q_{\delta,1}$ distribution, with
varying $\delta$.  We used 50 bootstrap replicates in each test, kept the
nominal size $\alpha = 10$, and calculated power by averaging over 25
independent graph pairs.  Despite the graphs having only 100 nodes,
Fig. \ref{fig:test} shows that our test has quite respectable power.

\paragraph{Comparison of Watts-Strogatz Graphs}

\begin{figure}[h]
    \includegraphics[width=25mm,height=25mm]{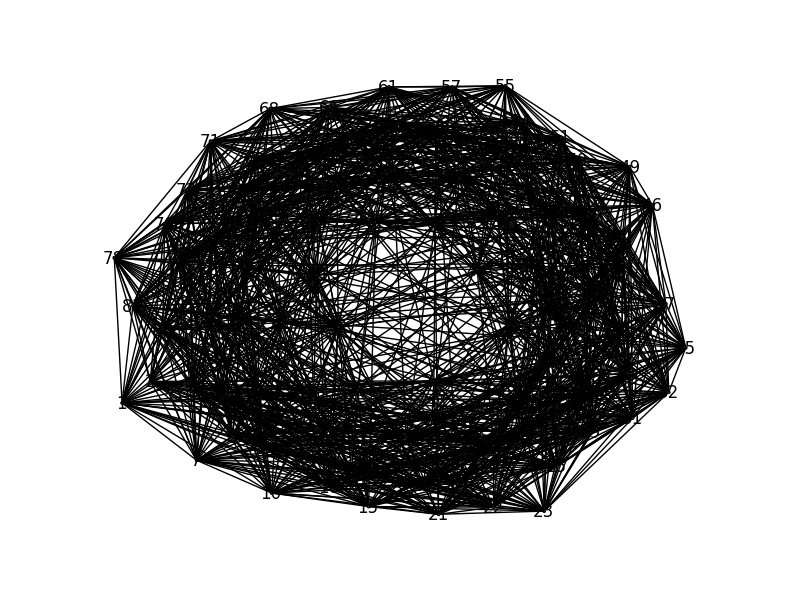}
    \includegraphics[width=25mm,height=25mm]{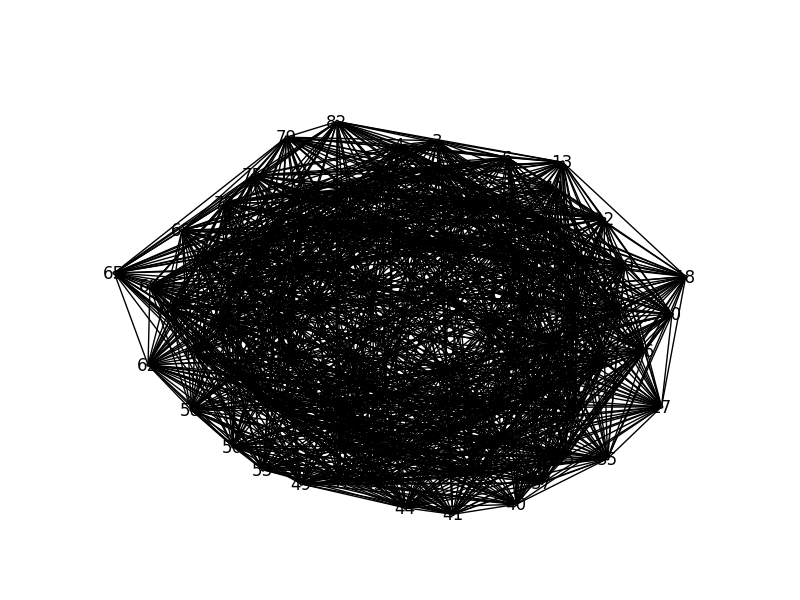}
    \includegraphics[width=25mm,height=25mm]{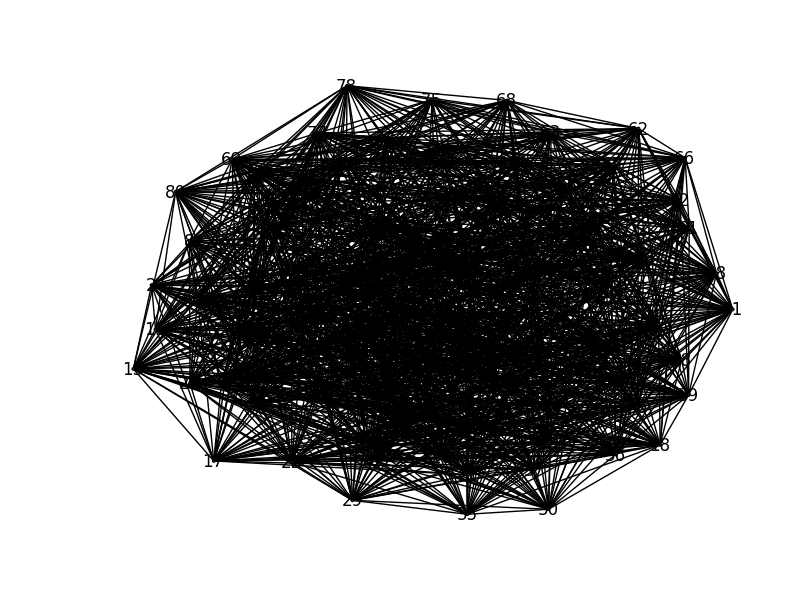}
    \includegraphics[width=25mm,height=25mm]{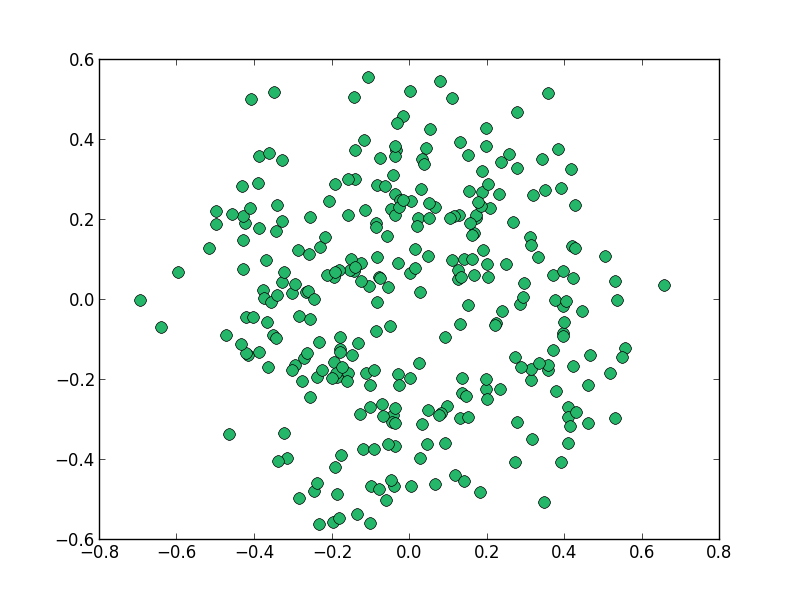}
    \includegraphics[width=25mm,height=25mm]{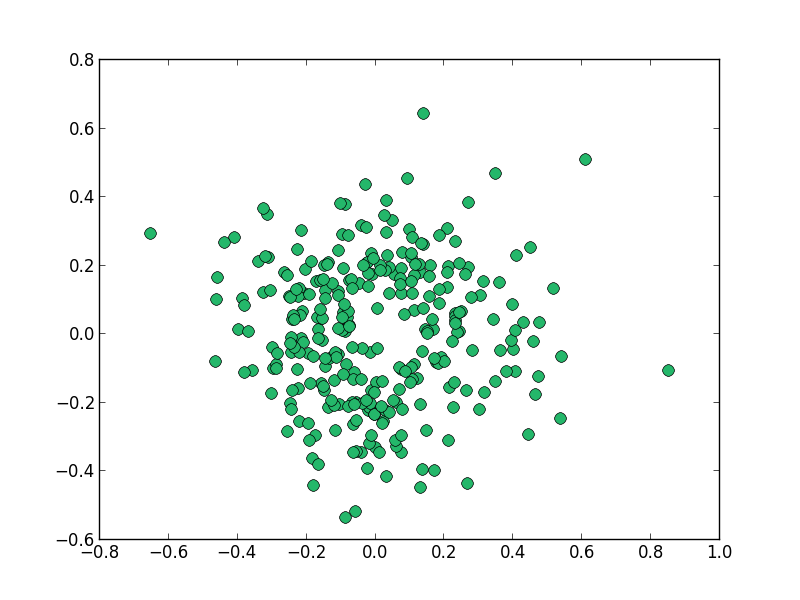}
    \includegraphics[width=25mm,height=25mm]{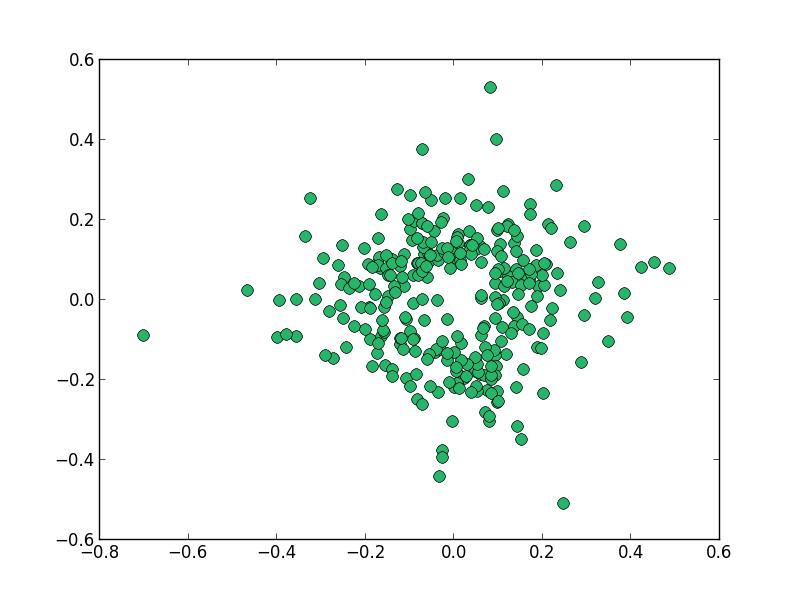}
\caption{{\bf Comparing Watts–-Strogatz models} Above, Watts-Strogatz graphs formed by re-wiring 1D ring lattices (85 nodes, branching factor 40) with probability $p$ per edge; from left to right $p=0.1, 0.2, 0.3$.  Below, embeddings of the graphs into $\hyperboloid{2}$.  At nominal $\alpha=0.1$, the power to detect these differences in $p$ was $1.0$ to within Monte Carlo error.}
  \label{fig:watts-strogatz}
\end{figure}

We have explained above, \S \ref{sec:motivation}, why we expect hyperbolic
latent space models to be reasonable ways of summarizing the structure of
complex networks.  However, they will also be more or less mis-specified for many
networks of interest.  We thus applied our methods to a class of graph
distributions which do {\em not} follow a hyperbolic latent space model, namely
Watts-Strogatz networks \citep{Watts-Strogatz-small-world}.  Our simulations
used 100 node networks, with the base topology being a 1D ring with a branching
factor of 40, and variable re-wiring probabilities.  These graphs show the
small-world property and high transitivity, but light-tailed degree
distributions.  Even in these cases, where the hyperbolic model is not the true
generator, our comparison method had almost perfect power
(Fig. \ref{fig:watts-strogatz}).

%% file: conclusion.tex
We have shown how nonparametric hyperbolic latent space models let us compare
the global structures of networks.  Our approach has its limits, and it may
work poorly when the networks being compared are very far from hyperbolic.
However, our experiments with Watts-Strogatz graphs show that it can detect
differences among graph distributions from outside our model class.  When we do
detect a change in structure, we have a model for each network, namely their
node densities, and the difference in node densities is an interpretable
summary of how the networks differ.  Many important directions for future work
are now open.  One is a better handling of sparse networks, perhaps through
some size-dependent modification of the link-probability function $W$, as in
\cite{Krioukov-et-al-hyperbolic-geometry}, or the sort of scaling of graphons
introduced in \cite{Borgs-Chayes-Zhao-sparse-graph-convergence}.  But this
should only extend our method's scope.

%% file: helgason.tex
The reader is referred to \citep{Terras-harmonic-1} for a general defintion of the Helgason-Fourier transform on a symmetric space.  
We specialize those constructions for $\hyperboloid{2}$, regarded in this section as the Poincar\'{e} half-plane with the metric
$$dz=\nicefrac{dx\;dy}{y^2}.$$ 

Let $f$ and $\phi$ denote smooth maps 
$$f:\hyperboloid{2}\rightarrow\mathbb{C},\quad\phi:\mathbb{R}\times\mathbb{SO}_2\rightarrow\mathbb{C}$$
with compact supports.

The \textit{Helgason-Fourier transform} $\mathcal{H}[f]$ is the map
$$\mathcal{H}[f]:\mathbb{C}\times\mathbb{SO}_2\rightarrow\mathbb{C}$$
where $\mathbb{SO}_2$ is the space of all rotation matrices $k_\theta$ of positive determinant determined by an angle $\theta$, sending each pair $(s,k_\theta)$ to the integral
\begin{equation*}
  \int_{\hyperboloid{2}}f(z)\mathrm{Im}(k_\theta(z))^{\bar{s}}\;dz
\end{equation*}

The \textit{inverse Helgason-Fourier transform} is the map
$$\mathcal{H}^{-1}[\phi]:\hyperboloid{2}\rightarrow\mathbb{C}$$
sending each $z$ to the integral
$$\int_{\R}\int_{\mathbb{SO}_2}\phi(\nicefrac{1}{2}+it,k_\theta)\mathrm{Im}(k_\theta(z))^{\nicefrac{1}{2}+it}\nicefrac{t\;\tanh(t)\;d\theta\;dt}{8\pi^2}$$

For each $f$, we have the identities 
\begin{equation}
\mathcal{H}^{-1}[\mathcal{H}[f]]=f,\quad\|f\|_2=\|\mathcal{H}[f]\|_2.
\label{plancherel}
\end{equation}

In a certain sense, the operation $\mathcal{H}$ takes convolutions to products in the following sense.
Let $g$ denote a compactly supported density on $\SL_2$ that is \textit{$\SO_2$-invariant} in the sense that $g(axb)=g(x)$ for all $a,b\in\SO_2$.
Then each $g$ passes to a well-defined density on $\hyperboloid{2}=\SL_2/\SO_2$, which, by abuse of notation, is also written as $g$.  
As a function on $\SL_2$, the convolution $g*f$ can be defined as the density on $\hyperboloid{2}$ defined by
$$(g*f)(z)=\int_{\SL_2}g(m)f(m^{-1}z)dm,$$
where the integral is taken with respect to the Haar measure on $\SL_2$, the measure on $\SL_2$ that is unique up to scaling and invariant under multiplication on the left or right by an element.
Then for all $f$ and $g$,
$$\mathcal{H}[g*f]=\mathcal{H}[g]\mathcal{H}[f].$$

Moreover, the Helgason-Fourier transform and its inverse each send real-valued functions to real-valued functions.
For details, the reader is referred to \citep{Terras-harmonic-1}.

%% file: computing.tex
We compute our test statistics as follows.  
Given a pair of graphs $G_1$ and $G_2$ (of possibly varying size), we use generalized multidimensional scaling to obtain coordinates for $G_1$ and $G_2$, functions
$$\phi_1:{V_1}\rightarrow\hyperboloid{2},\quad\phi_2:{V_2}\rightarrow\hyperboloid{2}$$
from the vertices $V_1$ of $G_1$ and $V_2$ of $G_2$.
In our power tests, we generate our graphs $G_1$ and $G_2$ by sampling $100$ points from two different densities on $\hyperboloid{2}$ and connecting those points according to the Heaviside step function, as outlined in Figure \ref{fig:generated.graphs}.
Very rarely, generalized multidimensional scaling fails in the sense that $\cosh$ applied to the distance matrix does not have two negative eigenvalues.  When such failure occurs during our power tests, we simply generate two new graphs of $100$ nodes from appropriate densities on $\hyperboloid{2}$ and attempt the embedding algorithm once again.  

We then let $f_i$ be the generalized kernel density estimator (\ref{eqn:H2.estimator}) on $\hyperboloid{2}$ determined by the points $\phi_i(V_i)$, where $n=\#\phi_i(V_i)$, $T=n^{-\nicefrac{1}{6}}$, and $h=\nicefrac{1}{n+100}$ (experience shows that for the $n$ values with which we are working, this choice of bandwidth $h$ works best.)
Thus the network models we estimate for $G_1$ and $G_2$ are the continuous latent space models $\hat{P}_1$ and $\hat{P}_2$ respectively determined by $f_1$ and $f_2$.  
The test statistic $d^*=d(\hat{P}_1,\hat{P}_2)$ we compute is
$$\|f_1-f_2\|_2=\int_{-T}^{+T}\int_{0}^{2\pi}(\mathcal{H}[f_2]-\mathcal{H}[f_1])^2 \nicefrac{t\tanh t}{8\pi^2}\;d\theta\;dt.$$
approximated by averaging the integrand over $100$ uniformly chosen pairs $(t,\theta)\in[-T,+T]\times[0,2\pi)$.
The integrand itself is computed as follows.
The Helgason-Fourier transform $\mathcal{H}[f_i]$ of (\ref{eqn:H2.estimator}) is the function
$$\mathcal{H}[f_i]:\mathbb{R}\times\mathbb{SO}_2\rightarrow\mathbb{R}$$
defined by the rule
\begin{equation}
  \label{eqn:integrand}
  \mathcal{H}[f_i](t,\theta)=\frac{1}{n}\sum_{i=1}^n\mathrm{Im}(k_\theta(Z_i))^{\nicefrac{1}{4}-t^2}e^{\nicefrac{1}{4}-t^2}
\end{equation}